\documentclass[hidelinks]{siamart220329}

\usepackage{lipsum}
\usepackage{amsfonts}
\usepackage{graphicx}
\usepackage{epstopdf}
\usepackage{enumerate}
\usepackage{hyperref}

\ifpdf
  \DeclareGraphicsExtensions{.eps,.pdf,.png,.jpg}
\else
  \DeclareGraphicsExtensions{.eps}
\fi

\newsiamremark{remark}{Remark}
\newsiamremark{hypothesis}{Hypothesis}
\crefname{hypothesis}{Hypothesis}{Hypotheses}
\newsiamthm{claim}{Claim}
\newsiamthm{example}{Example}
\newsiamthm{assumption}{Assumption}
\headers{Exact Convergence rate of Polyak step size}{M. Zamani, F. Glineur}

\title{Exact Convergence rate of the subgradient method by using Polyak step size
}

\author{Moslem Zamani\thanks{ICTEAM/INMA, Université catholique de Louvain, B-1348 Louvain-la-Neuve, Belgium
  (\email{moslem.zamani@uclouvain.be}).}
\and Fran\c{c}ois Glineur \thanks{ICTEAM/INMA \& CORE, Université catholique de Louvain, B-1348 Louvain-la-Neuve, Belgium
  (\email{Francois.Glineur@uclouvain.be}).
  }}
\usepackage{amsopn}

\DeclareMathOperator{\supp}{supp}
\usepackage{algorithm,algpseudocode}
\algnewcommand{\Inputs}[1]{%
  \State \textbf{Inputs:}
  \Statex \hspace*{\algorithmicindent}\parbox[t]{.8\linewidth}{\raggedright #1}
}
\algnewcommand{\Initialize}[1]{%
  \State \textbf{Initialization:}
  \Statex \hspace*{\algorithmicindent}\parbox[t]{.8\linewidth}{\raggedright #1}
}

\ifpdf
\hypersetup{
  pdftitle={An Example Article},
  pdfauthor={D. Doe, P. T. Frank, and J. E. Smith}
}
\fi

\begin{document}

\maketitle

\begin{abstract}
This paper studies the last iterate of subgradient method with Polyak step size  when applied to the minimization of a nonsmooth convex function with bounded subgradients. We show that the subgradient method with Polyak step size achieves a convergence rate $\mathcal{O}\left(\tfrac{1}{\sqrt[4]{N}}\right)$ in terms of the final iterate.  An example is provided to show that this rate is exact and cannot be improved. We introduce an adaptive Polyak step size for which the subgradient method enjoys a convergence rate $\mathcal{O}\left(\tfrac{1}{\sqrt{N}}\right)$ for the last iterate. Its convergence rate  matches exactly the lower bound on the performance of any black-box method on the considered problem class. Additionally, we propose an adaptive Polyak method with a momentum term, where the step sizes are independent of the number of iterates. We establish that the algorithm also attains the optimal convergence rate.  
We  investigate the alternating projection method. We derive a convergence rate $\left( \frac{2N }{ 2N+1 } \right)^N\tfrac{R}{\sqrt{2N+1}}$ for the last iterate,   where $R$ is a bound on the distance between the initial iterate and a solution. An example is also provided to illustrate the exactness of the rate.
\end{abstract}

\begin{keywords}
convex optimization, nonsmooth optimization, subgradient method, Polyak step size, convergence rate, last iterate, alternating projection method
\end{keywords}

\begin{MSCcodes}
90C25, 90C60, 49J52
\end{MSCcodes}


\section{Introduction}\label{intro}
The subgradient method is an iterative technique designed to solve nonsmooth convex optimization problems. Originally studied by Shor and others in the 1960s, the method is both straightforward and widely employed, continuing to be a topic of active research. Recent advancements have led to the development of new, more efficient variants capable of tackling a broader class of optimization problems, as discussed in \cite{davis2019stochastic, grimmer2023some} and related references.

Let \( X \subseteq \mathbb{R}^n \) be a closed convex set and let $f$ be a convex function whose domain, denoted by $\mathrm{dom} f$, includes \( X \). Consider the convex optimization problem:
\begin{align}\label{P}
\min_{x \in X} f(x).
\end{align}

At any point \( x \), the set of subgradients of the function \( f \), denoted by \( \partial f(x) \), is defined as:
$$
\partial f(x) = \{ g \in \mathbb{R}^n \mid f(y) \geq f(x) + \langle g, y - x \rangle \text{ for all } y \in \mathrm{dom} f \}.
$$
Throughout the paper, we assume that  the set of subdifferential $\partial f(x)$ is nonempty for every $x \in X$ and it is bounded, i.e. for some $B>0$, 
$$
\|g\|\leq B, \ \ \ \ \forall g\in\partial f(x), x\in X. 
$$
Moreover, we assume that the optimal solution set of problem \eqref{P}, $X^\star$, is non-empty. 
Algorithm \ref{SM} outlines the  projected subgradient method, where \( \Pi_X(\cdot) \) represents the Euclidean projection onto the set \( X \). Implementing this method typically involves defining a sequence of \( N \) step size, \( \{ h_k \}_{1 \le k \le N} \), where \( N \) denotes the number of iterations.

\begin{algorithm}
\caption{Projected subgradient method with generic step size}
\begin{algorithmic}
\smallskip
\State \textbf{Parameters:} number of iterations $N$.
\smallskip
\State \textbf{Inputs:} convex set $X$, convex function $f$ defined on $X$, initial iterate $x^1 \in X$.
\smallskip
\State For $k=1, 2, \ldots, N$ perform the following steps:\\
\begin{enumerate}
\item
Select a subgradient $g^k\in \partial f(x^k)$.
\item
Compute $x^{k+1}=\Pi_X\left(x^k-h_kg^{k}\right)$ using the step size $h_k$.
\end{enumerate}
\smallskip
\State \textbf{Output:} last iterate $x_{N+1}$ 
\end{algorithmic}
\label{SM}
\end{algorithm}

Employing appropriate step sizes is crucial for the effectiveness of the subgradient method. Indeed, it is directly influencing the convergence and  efficiency of the algorithm. If the step size is too large, the algorithm may overshoot the optimal point, leading to oscillations. Conversely, if the step size is too small, the algorithm may converge very slowly \cite{boyd}. Among the various strategies for choosing step size, the Polyak step size stands out due to its potential to significantly enhance convergence rates \cite{barre2020complexity, hazan2019revisiting}. The Polyak step size offers a significant advantage by dynamically adjusting the step size based on the current state of the optimization problem in question. 

 The Polyak step size is defined as,
\begin{align}\label{Polyak_t}
h_k = \frac{ t_k\left(f(x^k) - f^\star\right)}{\| g^k \|^2},
\end{align}
where \( f^\star \) denotes the optimal value of problem \eqref{P} and $t_k\in(0,2)$. Commonly,  $t_k=1$ is called the Polyak step size in the literature \cite{boyd, hazan2019revisiting}.  In the rest, we use  the Polyak step size for $t_k=1$. By adjusting the step size in accordance with the distance from the optimal value, the Polyak step size can lead to faster convergence compared to fixed or diminishing step size.

A notable limitation of the Polyak step size is the requirement to know \( f^\star \) in advance. For some class of problems including the convex feasibility problem or over-parameterized machine learning models, we may have access to the optimal value. Nevertheless, in many practical optimization problems, the exact optimal value is unavailable, which makes a direct application of the Polyak step size challenging. When \( f^\star \) is not available, one common approach is to estimate \( f^\star \). Despite this drawback, the Polyak step size remains a powerful tool when the optimal value can be reasonably estimated or is known from prior information \cite{bertsekas2015convex, hazan2019revisiting}. 
 
By employing step size \eqref{Polyak_t}, the subgradient method generates a Fej\'{e}r monotone sequence $\{x^k\}$ with respect to $X^\star$, that is,
$$
\|x^{k+1}-x^\star\|\leq \|x^{k}-x^\star\|, \ \ \ \forall x^\star\in X^\star, 
$$
see \cite[Theorem 1]{polyak1969minimization}. 
Polyak also \cite[Theorem 1]{polyak1969minimization} proved that if $\epsilon\leq t_k \leq 2-\epsilon$ for some $\epsilon \in (0,1)$, then
$$
\lim_{k\to\infty} f(x^k)=f^\star.
$$
The most natural question concerning this limit is its non-asymptotic  convergence rate. Additionally, this question is practically significant, as the final iterate is often regarded as the output of the subgradient algorithm \cite{jain2019making, last_iterate_sgd}. This subject is addressed for smooth problems, see \cite{barre2021worst, barre2020complexity, huang2024analytic} and references therein, but it is not investigated for non-smooth problems to the best knowledge of authors. Analyzing the final iterate of the subgradinet method is not straightforward. 
This difficulty is due to the fact that the subgradient method is not a descent method, that is, \( f(x_{k+1}) \le f(x_k) \) does not necessarily hold at each iterate. Therefore, most convergence rates are given in terms of the best iterate or an average of the iterates; see e.g. \cite{beck2017first, bertsekas2015convex, boyd}.

In this paper, we investigate this question and its related topics. In Section \ref{exact_Polyak}, we establish exact convergence rate the subgradient method with the step size $h_k = \frac{ f(x^k) - f^\star}{\| g^k \|^2}$. We show that the rate of convergence is $\mathcal{O}(\tfrac{1}{\sqrt[4]{N}})$. Furthermore, we provide an example demonstrating that the given bound cannot be improved.  In Section \ref{Sec_adap}, we present an adaptive Polyak step size for which the last-iterate convergence rate matches the established lower bound for black-box nonsmooth convex optimization problems, achieving \( f(x^N) - f(x^*) \le \frac{BR}{\sqrt{N+1}} \). However, the algorithm is not universal, that is, the step sizes are dependent on the number of iterations. We introduce an adaptive Polyak method with a momentum term that is universal and attains the optimal convergence rate.  

 Section \ref{Sec_alt_P} is devoted to convex feasibility problems. We develop algorithm introduced in Section \ref{Sec_adap} to this class of problems. We establish these projected based algorithms attain the optimal convergence rate for some class of problems. Moreover,  we study the alternating projection method. We derive a convergence rate $\left( \frac{2N }{ 2N+1 } \right)^N\tfrac{R}{\sqrt{2N+1}}$ for the last iterate,   where $R$ is a bound on the distance between the initial iterate and a solution. We give an example that demonstrates the exactness of the bound.

We conclude the section by recalling some notions and  preliminary results, which are used to derive results.

\subsection{Notations and  preliminaries}
We denote the $n$-dimensional Euclidean space by $\mathbb{R}^n$. We  denote
the Euclidean inner product and norm by $\langle \cdot, \cdot \rangle$ and $\| \cdot \|$, respectively.  We use $e$ to denote the vector of ones. In addition, $\{e_1, \dots, e_n\}$ represents the standard basis of $\mathbb{R}^n$. $A_{ij}$ denotes 
$(i, j)$-th entry of matrix $A$, and $A^T$ stands for the transpose of $A$. 
A symmetric matrix \( A \in \mathbb{R}^{n \times n} \) is called a Stieltjes matrix if we have for $i\neq j\in\{1, \dots, n\}$,
$$
A_{i,i}\geq 0, \ \ A_{ij}\leq 0, \ \ A^{-1}\geq 0,
$$
where the  matrix inequality is understood entrywise. Note that Stieltjes matrices are positive definite; see \cite[Theorem 2.3]{berman1994nonnegative}. The following lemma plays a key role in deriving the rates in this paper.

\begin{lemma}\label{Lemma1}\cite{zamani2023exact}
Let  $ x^\star\in X^\star$, $h_{N+1}>0$ and let $0<v_0\leq v_1\leq \dots \leq v_N\leq v_{N+1}$. If  Algorithm \ref{SM} with the starting point $x^1 \in X$  generates $\{(x^k, g^k)\}$, then
\begin{align}\label{UR_SM}
   & \sum_{k=1}^{N+1}  \left(h_kv_{k}v_{k-1}-(v_k-v_{k-1})\sum_{i=k+1}^{N+1}  h_iv_{i} \right) f(x^{k}) -
v_0 \sum_{k=1}^{N+1}  h_kv_{k}  f^\star 
 \leq \\
 \nonumber & \ \ \ \ 
  \tfrac{v_{0}^2}{2}\left \| x^{1}-x^\star\right\|^2 +\tfrac{1}{2}\sum_{k=1}^{N+1}  h_k^2v_{k}^2 \left \|g^{k}\right\|^2.
\end{align}
 \end{lemma}

The authors use Lemma \ref{Lemma1} to derive the exact convergence of subgradient method in terms of the last iterate for constant step size and length; see \cite{zamani2023exact}. For the extension of the given lemma to other settings, interested reader may refer to \cite{cai2024last, liu2024last, liu2023revisiting}. 

Due to the subgradient inequality, we get the following inequality for Algorithm \ref{SM},
\begin{align*}
f(x^{N+1})-f(x^\star)&=\sum_{k=1}^N \left( f(x^{k+1}-f(x^k) \right)+f(x^1)-f(x^\star)\leq \\
&
\sum_{k=1}^N \langle g^{k+1}, x^{k+1}-x^k \rangle+\langle g^1, x^1-x^\star \rangle.
\end{align*}
 One can infer from the inequality  for the worst-case function of Algorithm \ref{SM} equipped with step size \eqref{Polyak_t}
in the dimension-free analysis, we have $\|g^k\|=B$, $k\in\{1, \dots, N+1\}$.  We use this point in our analysis.
\section{Exact convergence rate of the Polyak step size ($t_k=1)$}\label{exact_Polyak}
In this section, we study the exact convergence rate of subgradient method equipped with $h_k = \frac{ f(x^k) - f^\star}{\| g^k \|^2}$ in terms of last iterate.  Before we get to the theorem, we need to present some lemmas. The following lemma is resulted from the Wallis product.

\begin{lemma}\label{Lemma_Polyak_0}
Let $1\leq k <N$. Then 
\begin{align}\label{Wallis}
    1< \tfrac{2k+1}{2k}\prod_{i=k+1}^N \tfrac{4i^2-1}{4i^2}< 2.
\end{align}
 \end{lemma}
\begin{proof}
Let $a_j=\int_0^{\tfrac{\pi}{2}} (\sin t)^jdt$. By using the Wallis product see e.g. \cite{lang2012first}, we have
$a_{2N}=\tfrac{\pi}{2}\prod_{i=1}^N \tfrac{2i-1}{2i}$, $a_{2N+1}=\prod_{i=1}^N \tfrac{2i}{2i+1}$ and 
$$
1\leq \frac{a_{2N}}{a_{2N+1}}\leq 1+\tfrac{1}{2N}.
$$
As $\prod_{i=k+1}^N \tfrac{4i^2-1}{4i^2}=\tfrac{a_{2N}/a_{2N+1}}{a_{2i}/a_{2i+1}}$ and $a_j$ is strictly decreasing, we get
$$
\tfrac{2k+1}{2k}\prod_{i=k+1}^N \tfrac{4i^2-1}{4i^2}\leq (\tfrac{2k+1}{2k})( 1+\tfrac{1}{2N})< 2.
$$
On the other hand, 
$$
\tfrac{2k+1}{2k}\prod_{i=k+1}^N \tfrac{4i^2-1}{4i^2}> (\tfrac{2k+1}{2k})(1+\tfrac{1}{2k})^{-1}> 1.
$$
\end{proof}
 \begin{lemma}\label{Lemma_PD_Polyak}
Let $N$ be given and $a_k$ defined as follows,
$$
a_k=\prod_{i=N+1-k}^N  \left(\tfrac{2i+1}{2i}\right)\left(\tfrac{4i^2}{4i^2-1}\right)^{N+1-i-k}, \ \ \ \ \ k\in \{1, \dots, N\}.
$$
and $a_0=1$. Suppose that $Q$ is an $N\times N$ symmetric matrix given as, 
 $$
Q_{ij}=\begin{cases}
 2a_i a_{i-1}-a_i^2  & i=j\\
 \left(a_{\min(i,j)-1}-a_{\min(i,j)}\right)a_{\max({i,j})}   & i\neq j.
\end{cases}
 $$
 Then matrix $Q$ is positive definite.
 \end{lemma}
 \begin{proof}
To establish the semi-definiteness of $Q$, we will show that the matrix is a Stieltjes matrix.  By the definition of $a_k$, we have 
$$
\tfrac{a_k}{a_{k-1}}=\tfrac{2(N+1-k)+1}{2(N+1-k)}\prod_{i=N+2-k}^N \left(\tfrac{4i^2-1}{4i^2}\right).
$$
Hence, Lemma \ref{Lemma_Polyak_0} implies that the diagonal elements of $Q$ are positive and the off-diagonal elements are negative. 
We define $p_k$ and $q_k$ for $k\in\{1, \dots, N\}$ as follows,
$$
p_k=\prod_{i=N+1-k}^N \tfrac{2i+1}{2i}, \ \ \ 
q_k=\prod_{i=N+2-k}^N \tfrac{2i}{2i-1},
$$
and $q_1=1$. Since $p_k=2(N+1-k)p_k-2(N-k)p_{k+1}$ for $k\in\{1, \dots, N-1\}$, we have 
\begin{align}\label{pk}
&\sum_{i=1}^{k-1} p_i=2N+1-2(N+1-k)p_k, \\
\nonumber
&\sum_{i=k+1}^N p_i=\sum_{i=k+1}^{N-1} p_i+p_N=2(N-k)p_{k+1}-p_N=(2(N-k)+1)p_k-p_N.
\end{align}
Analogously, we have 
\begin{align}\label{qk}
\sum_{i=1}^{k-1} q_i=1+\sum_{i=2}^{k-1} q_i=1+2N-(2(N+2-k)-1)q_k.
\end{align}
 Suppose that $y\in\mathbb{R}^N_+$ is given by 
\begin{align*}
    y_k=\prod_{i=N+1-k}^N \left(\tfrac{4i^2}{4i^2-1}\right)^{i+k-N-1} , \ \ \ k\in\{1, \dots, N\}.
\end{align*}
It is seen that $p_k=a_ky_k$, $q_k=a_{k-1}y_k$ and $a_{k-1}=\tfrac{a_k q_k}{p_k}$. By employing \eqref{pk} and \eqref{qk}, we get
\begin{align}\label{opt_y}
\sum_{j=1}^N & Q_{kj}y_j=\sum_{j=1}^{k-1} Q_{kj}y_j+ Q_{kk}y_k+\sum_{j=k+1}^N Q_{kj}y_j=
a_k\sum_{j=1}^{k-1} (a_{j-1}-a_j)y_j+ \\
\nonumber & (2a_{k-1}-a_k)p_k+(a_{k-1}-a_k)\sum_{j=k+1}^N  a_jy_j=(a_k-\tfrac{a_kq_k}{p_k})p_N=(a_k-a_{k-1})p_N.
\end{align}
So $Qy>0$, and we can infer from \cite[Theorem 2.3 (I28)]{berman1994nonnegative} that $Q$ is a Stieltjes matrix. Consequently, $Q$ is positive definite, completing the proof.
 \end{proof}
 
\begin{theorem}\label{T_Polyak}
Let $f$ be a convex function with $B$-bounded subgradients on a convex set $X$. Consider $N$ iterations of Algorithm \ref{SM} with $h_k = \frac{ f(x^k) - f^\star}{\| g^k \|^2}$ starting from an initial iterate $x^1 \in X$ satisfying $\|x^1 - x^\star\| \le R$ for some minimizer $x^\star$. 
We have that the last iterate $x_{N+1}$ satisfies
\begin{align}\label{R_Polyak}
f(x^{N+1})-f(x^\star)\leq \tfrac{BR}{\sqrt{2N+1}}\prod_{i=1}^N \left(\tfrac{4i^2}{4i^2-1}\right)^{i}.
\end{align}
\end{theorem}
\begin{proof}
Without loss of generality, we assume $f^\star=0$. For the convenience of notation, let $f^k=f(x^k)$ for $k\in\{1, \dots, N\}$. As discussed in Introduction for the worst-case scenario, we have $\|g^k\| = B$ for $k \in \{1, \dots, N+1 \}$. Hence, we assume that $\|g^k\| = B$ in the argument.
Let $v_0=\tfrac{\sqrt{B}}{\sqrt[4]{R^2(2N+1)}}\prod_{i=1}^N \left(\tfrac{4i^2}{4i^2-1}\right)^{\tfrac{i}{2}}$ and $v_k$ be given as follows,
$$
v_k=v_0\prod_{i=N+1-k}^N  \left(\tfrac{2i+1}{2i}\right)\left(\tfrac{4i^2}{4i^2-1}\right)^{N+1-i-k}, \ \ \ \ \ k\in \{1, \dots, N\},
$$
and $v_{N+1}=v_N$. It is readily seen that $0<v_0\leq v_1\leq \dots \leq v_{N}\leq v_{N+1}$. Suppose that $h_{N+1}=\tfrac{1}{v_N^2}$. 
By  Lemma \ref{Lemma1}, we have
\begin{align*}
& f(x^{N+1})+\sum_{k=1}^{N}  \left(\tfrac{v_k v_{k-1}}{B^2}f^k-\tfrac{v_k-v_{k-1}}{ B^2 }\sum_{i=k+1}^{N}  v_if^i -\tfrac{ v_k-v_{k-1}}{v_N} \right) f^k 
 \leq \\
 \nonumber & \leq  \tfrac{1}{2B^2}\sum_{k=1}^{N} (v_kf^k)^2+\tfrac{B^2}{2v_N^2}+\tfrac{v_0^2}{2}\left\|x^1-x^\star\right\|^2.
\end{align*}
 As $\|x^1-x^\star\|\leq R$, the inequality may be written as
 \begin{align*}
    f(x^{N+1})\leq -\tfrac{v_0^2}{2B^2} F^TQF+b^TF+\tfrac{B^2}{2v_N^2}+\tfrac{v_0^2R^2}{2}:=q(F),
 \end{align*}
 where $F=\begin{pmatrix}
     f^1 & \dots & f^N
 \end{pmatrix}^T$,
 the symmetric matrix $Q$ given in Lemma \ref{Lemma_PD_Polyak}
 and 
 $$
b_k=\tfrac{ v_k-v_{k-1}}{v_N}, \ \ \ k\in\{1, \dots, N\}.
 $$
 By Lemma \ref{Lemma_PD}, $Q$ positive definite. So, $q$ is a concave function. To establish bound \eqref{R_SMO}, it suffices to determine the optimal value of the following problem
 $$
 \max q(F).
 $$
For $p_N$ given in Lemma \ref{Lemma_PD_Polyak}, we have $p_N=\tfrac{B\sqrt{2N+1}}{v_N v_0 R}$. Hence,   by \ref{opt_y}, we get $\nabla q(\bar F)=0$, where 
 \begin{align*}
    \bar F_k=\tfrac{BR}{\sqrt{2N+1}}\prod_{i=N+1-k}^N \left(\tfrac{4i^2}{4i^2-1}\right)^{i+k-N-1} , \ \ \ k\in\{1, \dots, N\}.
\end{align*}
On the other hand, we have
\begin{align}\label{FFFF}
\nonumber v_N^2=&
\tfrac{B}{R\sqrt{2N+1}}\prod_{i=1}^N  \left(\tfrac{2i+1}{2i}\right)^2\left(\tfrac{4i^2}{4i^2-1}\right)^{2-i}=
\tfrac{B}{R\sqrt{2N+1}}\prod_{i=1}^N  \left(\tfrac{2i+1}{2i-1}\right)\left(\tfrac{4i^2}{4i^2-1}\right)^{1-i}\\
&=
\tfrac{B\sqrt{2N+1}}{R}\prod_{i=1}^N \left(\tfrac{4i^2}{4i^2-1}\right)^{1-i}.
\end{align}
By \eqref{pk}, \eqref{qk}, \eqref{FFFF} and  $\nabla q(\bar F)=0$, we get
  \begin{align*}
  q(\bar F)= \tfrac{1}{2}b^T\bar F+\tfrac{B^2}{2v_N^2}+\tfrac{v_0^2R^2}{2}=\tfrac{v_0^2R^2}{2}-\tfrac{B^2}{2v_N^2}+\tfrac{B^2}{2v_N^2}+\tfrac{v_0^2R^2}{2}=\tfrac{BR}{\sqrt{2N+1}}\prod_{i=1}^N \left(\tfrac{4i^2}{4i^2-1}\right)^{i},
 \end{align*}
  which completes the proof. 
\end{proof}

In what follows, we construct an example that establishes the given bound in Theorem \ref{T_Polyak} is tight. Indeed, we introduce an example for which \eqref{R_Polyak} holds as equality. Before we get to the theorem we need to present some lemmas.  

 \begin{lemma}\label{Lemma_L1}
 Let $a_{N}=\prod_{i=1}^N\left( \tfrac{4i^2}{4i^2-1} \right)^{i}$,
 then $a_{N+1}=\Theta(\sqrt[4]{N})$.
 \end{lemma}
 \begin{proof}
Assume that $i\in\mathbb{N}$. We have 
$$
\tfrac{4i^2}{4i^2-1}=1+\tfrac{1}{4i^2-1}\geq 1+\tfrac{1}{4i^2}.
$$
As $\ln(1+x)\geq x-\tfrac{1}{2}x^2$ for $x\geq 0$, we get 
\begin{align*}
\ln(a_N)=&\sum_{i=1}^N i\ln( \tfrac{4i^2}{4i^2-1})\geq \sum_{i=1}^N i\ln( 1+\tfrac{1}{4i^2})\geq \sum_{i=1}^N (\tfrac{1}{4i}-\tfrac{1}{16i^3})\geq 
 \tfrac{1}{4}\ln(N-1)-\tfrac{1}{8}.
\end{align*}
Hence, $a_{N+1}=\Omega(\sqrt[4]{N})$. On the other hand, we have $\tfrac{1}{1-x}\leq 1+x+(\tfrac{4}{3})^3 x^2$ for $x\in [0, \tfrac{1}{4}]$. So,
$$
\ln(\tfrac{4i^2}{4i^2-1})\leq \ln\left(1+\tfrac{1}{4i^2}+\tfrac{16}{27i^4}\right)\leq \tfrac{1}{4i^2}+\tfrac{16}{27i^4}.
$$
Analogous to the former case, we get $a_{N+1}=\mathcal{O}(\sqrt[4]{N})$
and the proof is complete. 
 \end{proof}

By Theorem \ref{T_Polyak} and Lemma \ref{Lemma_L1}, one can infer that the subgarident method with the Polyak step size has a convergence rate $\mathcal{O}\left(\tfrac{1}{\sqrt[4]{N}}\right)$.

 \begin{lemma}\label{Lemma_L2}
 Let $\{a_k\}$ be given as follows, 
 \begin{align*}
    a_k=\prod_{i=N+1-k}^N \tfrac{4i^2}{4i^2-1} , \ \ \ k\in\{1, \dots, N\},
\end{align*}
and let $A=\begin{pmatrix}
    1 & ce^T\\
    ce & Q
\end{pmatrix}$, where $Q\in \mathbb{R}^{(N+1)\times (N+1)}$ given by
 $$
Q_{ij}=\begin{cases}
1   & i=j\\
1-a_{\min(i,j)}  & i\neq j,
\end{cases}
 $$
 and $c=\tfrac{1}{\sqrt{2N+1}}$.
 Then  there exists $(N+1)\times(N+2)$ matrix $G$ with $A=G^TG$.
 \end{lemma}
 \begin{proof}
First we  prove the positive definiteness of $Q$. To this end, we show that $Q$ is a Stieltjes matrix. It readily seen that the off-diagonal components of $Q$ are negative. Let vector $y\in\mathbb{R}^{N+1}_+$ be given as follows,
$$
y_k=\prod_{i=N+1-k}^N \tfrac{2i+1}{2i} , \ \ \ k\in\{1, \dots, N\},
$$
and $y_{N+1}=y_N$. For $k\in\{1, \dots, N-1\}$, we have $y_k=2(N+1-k)y_k-2(N-k)y_{k+1}$. By using the telescopic sum, we get
\begin{align}\label{S_y}
    \sum_{i=k}^{N+1} y_k= 2\sum_{i=k}^{N-1}\left( (N+1-i)y_i-(N-i)y_{i+1}\right)+2y_N=2(N+1-k)y_k.
\end{align}
On the other hand, we have 
$$
y_ka_k=\prod_{i=N+1-k}^N \tfrac{2i}{2i-1} , \ \ \ k\in\{1, \dots, N\},
$$
and $y_ka_k=(2(N+1-k)-1)y_ka_k-(2(N-k)-1)y_{k+1}a_{k+1}$ for $k\in\{1, \dots, N-1\}$. For $k\in\{1, \dots, N-1\}$, we have
\begin{align}\label{S_ya}
    \sum_{i=1}^{k} y_k a_k= (2N-1)y_1 a_ 1-(2(N-k)-1)y_{k+1}a_{k+1}.
\end{align}
By using \eqref{S_y} and \eqref{S_ya}, we get
$$
\sum_{j=1}^{N+1} Q_{kj} y_j=\sum_{j=1}^{k-1} (1-a_j)y_j+y_k+(1-a_k)\sum_{j=k+1}^{N+1}  y_j=1.
$$
So $Qy=e$ and  \cite[Theorem 2.3 (I28)]{berman1994nonnegative} implies that $Q$ is a Stieltjes matrix. Consequently, $Q$ is positive definite and may be written as $Q=R^T R$ for some non-singular matrix $R\in\mathbb{R}^{(N+1)\times (N+1)}$. Since $R$ is non-singular, there exists $\bar x$ with $R^T \bar x=\tfrac{1}{\sqrt{2N+1}} e$. In addition, we have
$$
\| \bar x\|^2=\tfrac{1}{2N+1}e^T R^{-1}R^{-T}e=\tfrac{1}{2N+1}e^TQ^{-1}QQ^{-1}e=\tfrac{1}{2N+1}y^T e=1,
$$
where the last equality follows from \eqref{S_y}. For matrix $G=\begin{pmatrix}
    \bar x & R
\end{pmatrix}$,
we have  $A=G^TG$ and the proof is complete. 
 \end{proof}

The next theorem shows that the convergence rate given in Theorem \ref{T_Polyak} is exact. 

\begin{theorem}\label{Ex_Poly_E}
Let $N\in\mathbb{N}$. There exists a convex function $f:\mathbb{R}^{N+1}\to\mathbb{R}$ with $1$-bounded subgradients and a minimum point $x^\star$, for which  
\begin{align}
f(x^{N+1})-f(x^\star)=\tfrac{1}{\sqrt{2N+1}}\prod_{i=1}^N \left(\tfrac{4i^2}{4i^2-1}\right)^{i},
\end{align}
where $x^{N+1}$ is generated by Algorithm  \ref{SM} with $h_k = \frac{ f(x^k) - f^\star}{\| g^k \|^2}$, initial point $x^1$ and $\|x^1-x^\star\|= 1$.
\end{theorem}
\begin{proof}
Let $f^i$ be given as follows,
\begin{align*}
    f^k=\tfrac{1}{\sqrt{2N+1}}\prod_{i=N+1-k}^N \left(\tfrac{4i^2}{4i^2-1}\right)^{i+k-N-1} , \ \ \ k\in\{1, \dots, N+1\},
\end{align*}
and let $g^k$ denote $k+1$ column matrix $G$ given in Lemma \ref{Lemma_L2} and $z^1$ stands for the first column of $G$. It is seen that 
\begin{align}\label{f_a}
    f^{k+1}=a_kf^k,  k\in\{1, \dots, N\},
\end{align}
for $a_k$ given in Lemma \ref{Lemma_L2}. In addition, $\langle g^k, x^1\rangle=f^1$ for $k\in\{1, \dots, N+1\}$. We define 
$$
z^k=z^1-\sum_{i=1}^{k-1} f^i g^i, \ \ \ k\in\{1, \dots, N+1\}.
$$
Suppose that $g^{N+2}=z^{N+2}=0$ and $f^{N+2}=0$. Consider $f:\mathbb{R}^{N+1}\to\mathbb{R}$ given by $f(x)=\max_{1\leq k\leq N+2} f^k+\langle g^k, x-z^k \rangle$. By lemma \ref{Lemma_L2} and \eqref{f_a}, we have 
\begin{align*}
f^k+\langle g^k, 0-z^k \rangle&=f^k-\langle g^k, z^1 \rangle+\sum_{i=1}^{k-1} f^i\langle g^k, g^i \rangle=f^k-f^1+\sum_{i=1}^{k-1} f^i(1-a_i)\\
& =f^k-f^1+\sum_{i=1}^{k-1} (f^i-f^{i+1})=0.
\end{align*}
Hence,  $0\in \text{argmin} f(x)$ and $f(0)=0$. For $j<k$, we have
\begin{align}\label{f_k_j}
\nonumber f^k+\langle g^k, z^j-z^k \rangle&=f^k+\sum_{i=j}^{k-1} f^i\langle g^k, g^i \rangle=f^k+\sum_{i=j}^{k-1} f^i(1-a_i)\\
&= f^k+\sum_{i=j}^{k-1} (f^i-f^{i+1})=f^j.
\end{align}
Analogously, we obtain for $j>k$,
\begin{align}\label{f_j_k}
\nonumber f^k+\langle g^k, z^j-z^k \rangle&=f^k-\sum_{i=k}^{j-1} f^i\langle g^k, g^i \rangle=\sum_{i=k+1}^{j-1} (a_k-1)f^i\leq \sum_{i=k+1}^{j-1} (a_i-1)f^i\\
=f^j-f^{k+1}.
\end{align}
The inequality follows from the point that $a_k$ is increasing. By \eqref{f_j_k} and \eqref{f_k_j}, we have 
$$
f(z^k)=f^k, \ \ k\in\{1, \dots, N+2\}.
$$
Hence Algorithm  \ref{SM} with $h_k = \frac{ f(x^k) - f^\star}{\| g^k \|^2}$, initial point $x^1$ and selecting $g^k$ as a subgradient at iterate $k$ generates $x^{k+1}=z^{k+1}$, $k\in\{1, \dots, N\}$. Hence, $f(x^{N+1})=f^{N+1}$ and the proof is complete.
\end{proof}

One may wonder how we got the results in the section. In our previous work \cite[Example 3.5]{zamani2023exact}, we observed that the worst-case function for the subgradient with a fixed step size satisfies the following system
\begin{align*}
& f(x^i)=f(x^j)+\langle g^j, x^i-x^j\rangle, \ \ \ \ 1\leq i< j\leq N+1\\
& f(x^\star)=f(x^j)+\langle g^j, x^\star-x^j\rangle, \ \ \ \ 1\leq j\leq N+1.
\end{align*}
We solved the system for the Polyak step size and constructed the function given in Theorem \ref{Ex_Poly_E}. Moreover, we derived Theorem \ref{T_Polyak} by choosing suitable $v_k$'s and $h_{N+1}$. 
\section{Adaptive Polyak step size}\label{Sec_adap}
Recently, the authors \cite{zamani2023exact} introduced a new optimal subgradient method that achieves the best possible last-iterate accuracy. Strictly speaking, consider  Algorithm \ref{SM} with the initial point $x^1$ satisfying $\|x^1-x^\star\|\leq R$.  We proved that Algorithm \ref{SM} by employing 
\begin{align}\label{h_OPT}
    h_k=\tfrac{R(N+1-k)}{\|g^k\|\sqrt{(N+1)^3}}, \ \ k\in\{1, \dots, N\},
\end{align}
achieves a convergence rate $f(x^{N+1})-f(x^\star)\leq \tfrac{BR}{\sqrt{N+1}}$.   It is worth noting that we have  the following lower bound for the class of problems in question,
\begin{align}\label{Lower}
f(x^{N+1})-f(x^\star)\geq \frac{BR}{\sqrt{N+1}};
\end{align}
see \cite{drori2016optimal}. For more discussion on the step size \eqref{h_OPT} both theoretically and numerically,  we refer the reader to \cite{defazio2023and, defazio2024road}.

Inspired by \eqref{h_OPT}, we introduce Algorithm \ref{OSM}, an adaptive Polyak step size. Furthermore, we establish that this algorithm attains the best possible worst-case convergence rate.

\begin{algorithm}
\caption{adaptive Polyak step size}
\begin{algorithmic}
\smallskip
\State \textbf{Parameters:} number of iterations $N$
\smallskip
\State \textbf{Inputs:} convex set $X$, convex function $f$ defined on $X$ with $B$-bounded subgradients and the optimal value $f^\star$, initial iterate $x^1 \in X$ 
\smallskip
\State For $k=1, 2, \ldots, N$ perform the following steps:\\
\begin{enumerate}
\item
Select a subgradient $g^k\in \partial f(x^k)$.
\item
Set $x^{k+1}=\Pi_X\left(x^k- h_k  g^{k}\right)$ using step size $h_k=\tfrac{(N+1-k)(f(x^k)-f^\star)}{(N+1) \|g^k\|^2}$.
\end{enumerate}
\smallskip
\State \textbf{Output:} last iterate $x_{N+1}$ 
\end{algorithmic}
\label{OSM}
\end{algorithm}

As mentioned in Introduction the sequence $\{x^k\}$ generated by Algorithm \ref{OSM} is F\'ejer monotone w.r.t. $X^\star$; see \cite[Theorem 1]{polyak1969minimization}. Before we establish the convergence rate of Algorithm \ref{OSM}, we need to present some lemma.

 \begin{lemma}\label{Lemma_PD}
Let $a_k=\tfrac{1}{N+1-k}$ for $k\in\{0, \dots, N\}$ and $a_{N+1}=a_N$.
Suppose that $Q$ is an $(N+1)\times (N+1)$ symmetric matrix defined as follows, 
 $$
Q_{ij}=\begin{cases}
 2a_{i-1}-a_0   & i=j\\
 a_{\min(i,j)-1}-a_{\min(i,j)}   & i\neq j.
\end{cases}
 $$
 Then matrix $A=Q-e_{N+1}e_{N+1}^T$ is positive semi-definite.
 \end{lemma}
 \begin{proof}
To prove the positive semi-definiteness of $A$, we will demonstrate that the given matrix is diagonally dominant. For $k<N+1$, we have 
\begin{align*}
\sum_{\substack{j=1\\ j\neq k}}^{N+1} |A_{kj}|=&\sum_{j=1}^{k-1}(a_j-a_{j-1})+\sum_{j=k+1}^{N+1} (a_k-a_{k-1})=a_{k-1}-a_0+\\
&(N+1-k)(a_k-a_{k-1})=2a_{k-1}-a_0=A_{kk}.
\end{align*}
As columns $N$ and $N+1$ have same elements, we have $A_{N+1,N+1}= \sum_{j=1}^{N} |A_{N+1,j}|$
and the proof is complete. 
 \end{proof}

\begin{theorem}\label{T_OSM}
Let $f$ be a convex function with $B$-bounded subgradients on a convex set $X$. Consider $N$ iterations of Algorithm \ref{OSM} starting from an initial iterate $x^1 \in X$ satisfying $\|x^1 - x^\star\| \le R$ for some minimizer $x^\star$. 
We have that the last iterate $x_{N+1}$ satisfies
\begin{align}\label{R_SMO}
f(x^{N+1})-f(x^\star)\leq \frac{BR}{\sqrt{N+1}}.
\end{align}
\end{theorem}
\begin{proof}
 We may assume $f^\star=0$. For the convenience of notation, let $f^k=f(x^k)$ for $k\in\{1, \dots, N+1\}$. 
As for the worst-case problem, we have $\|g^k\| = B$ for $k \in \{1, \dots, N+1\}$. We proceed the argument with the assumption that $\|g^k\| = B$.
Let $v_k$'s be given as follows,
$$
v_k=\frac{\sqrt[4]{(N+1)^3}}{N+1-k}\sqrt{\frac{B}{R}}, \ \ \ \ \ k\in \{0, \dots, N\},
$$
and $v_{N+1}=v_N$. It is seen that $0<v_0\leq v_1\leq \dots \leq v_{N}\leq v_{N+1}$. Suppose that $h_{N+1}=\tfrac{f^{N+1}}{(N+1) B^2}$. 
By  Lemma \ref{Lemma1}, since $h_kv_k=\tfrac{v_0}{B^2}$ for $k\in\{1, \dots, N+1\}$, we get
\begin{align*}
& \sum_{k=1}^{N+1}  \left(\tfrac{v_0 v_{k-1}}{B^2}f^k-\tfrac{ v_0(v_k-v_{k-1}) }{ B^2 }\sum_{i=k+1}^{N+1}  f^i \right) f^k 
 \leq 
  \tfrac{v_0^2}{2B^2}\sum_{k=1}^{N+1} (f^k)^2+
\tfrac{B}{2R\sqrt{N+1}}\left\|x^1-x^\star\right\|^2.
\end{align*}
 As $\left\|x^1-x^\star\right\|\leq R$, the above inequality may be written as
 \begin{align*}
   \tfrac{\sqrt{N+1}}{2BR}(f^{N+1})^2+\tfrac{\sqrt{N+1}}{2BR} F^TAF\leq 
\tfrac{BR}{2\sqrt{N+1}},
 \end{align*}
 where $F=\begin{pmatrix}
     f^1 & \dots & f^{N+1}
 \end{pmatrix}^T$,
 the symmetric matrix $A$ given in Lemma \ref{Lemma_PD}.
Since $A$ is positive semi-definite, the above inequality implies the desired bound and the proof is complete. 
\end{proof}

As seen, the rate given in Theorem \ref{T_OSM} exactly matches lower bound \eqref{Lower}. Furthermore, the step sizes in Algorithm \ref{OSM} are dependent on the number of iterations, $N$. This raises the question of whether there exists a universal optimal algorithm whose convergence rate matches lower bound  \eqref{Lower}. Analogous to \cite[Section 5.3]{zamani2023exact}, one can show that it is not possible.  An alternative approach to achieving an optimal subgradient method in this context could be the incorporation of a momentum term. In what follows, we introduce Algorithm \ref{OSM+M} with the optimal  convergence rate. Note that Algorithm \ref{OSM+M} needs the Lipschitz constant of $f$ on $X$.

\begin{algorithm}
\caption{Optimal Polyak step size method with a momentum term}
\begin{algorithmic}
\smallskip
\State \textbf{Parameters:} number of iterations $N$
\smallskip
\State \textbf{Inputs:} convex set $X$, convex function $f$ defined on $X$ with $B$-bounded subgradients and the optimal value $f^\star$, initial iterate $x^1 \in X$ 
\smallskip
\State Set $x^0=x^1$. For $k=1, 2, \ldots, N$ perform the following steps:\\
\begin{enumerate}
\item
Select a subgradient $g^k\in \partial f(x^k)$.
\item
Set $x^{k+1}=\Pi_X\left(x^k- \tfrac{f(x^k)-f^\star}{(k+1) B^2}  g^{k}+\tfrac{k-1}{k+1}\left(x^k-x^{k-1}\right)\right)$.
\end{enumerate}
\smallskip
\State \textbf{Output:} last iterate $x_{N+1}$ 
\end{algorithmic}
\label{OSM+M}
\end{algorithm}

\begin{theorem}\label{T_OSM+M}
Let $f$ be a convex function with $B$-bounded subgradients on a convex set $X$. Consider $N$ iterations of Algorithm \ref{OSM+M} starting from an initial iterate $x^1 \in X$ satisfying $\|x^1 - x^\star\| \le R$ for some minimizer $x^\star$. 
We have that the last iterate $x_{N+1}$ satisfies
\begin{align*}
f(x^{N+1})-f(x^\star)\leq \frac{BR}{\sqrt{N+1}}.
\end{align*}
\end{theorem}
\begin{proof}
 Suppose that $f^k=f(x^k)$ for $k\in\{1, \dots, N+1\}$ and, without loss of generality,  we assume that $f^\star=0$. By subgradient inequality, we have for $k\in\{1, \dots, N+1\}$,
\begin{align}\label{F1+}
&(k-1)\left( f^{k-1}-f^k-\langle g^k, x^{k-1}-x^k\rangle \right)+\left( -f^k-\langle g^k, x^\star-x^k\rangle \right)=\\
\nonumber & (k-1)f^{k-1}-kf^k+\langle g^k, kx^k-(k-1)x^{k-1}- x^\star\rangle\geq 0.
\end{align}
By multiplying \eqref{F1+} by $\tfrac{f^k}{B^2}$, we get for $k\in\{1, \dots, N\}$,
\begin{align*}
& \tfrac{(k-1)f^{k-1}f^k}{B^2}-\tfrac{k(f^k)^2}{B^2}\geq \langle \tfrac{-f^k}{B^2}g^k, kx^k-(k-1)x^{k-1}- x^\star\rangle=
\\
& -\tfrac{1}{2}\left\| kx^k-(k-1)x^{k-1}- x^\star  \right\|^2
-\tfrac{(f^k)^2}{2B^2}+\tfrac{1}{2}(k+1)^2\\
 &\left\| x^k- \tfrac{f^k}{(k+1) B^2}  g^{k}+\tfrac{k-1}{k+1}\left(x^k-x^{k-1}\right)-\left(\tfrac{k}{k+1}x^k+\tfrac{1}{k+1}x^\star  \right) \right\|^2.
\end{align*}
Due to the non-expansive property of the projection operator, we obtain
\begin{align*}
  {2(k-1)f^{k-1}f^k}-&{(2k-1)(f^k)^2}\geq  B^2\left\| (k+1)x^{k+1}-kx^k- x^\star  \right\|^2-\\
  & B^2\left\| kx^k-(k-1)x^{k-1}- x^\star  \right\|^2.
\end{align*}
Summing above inequalities for $k\in\{1, \dots, N\}$, we get 
\begin{align}\label{F2+}
\nonumber -N(f^N)^2-\sum_{k=2}^N (k-1)\left(  f^k-f^{k-1} \right)^2&\geq
 B^2\left\| (N+1)x^{N+1}-Nx^N- x^\star  \right\|^2-
\\ &
B^2\left\| x^1-x^\star  \right\|^2.
\end{align}
Let $g^{N+1}\in\partial f(x^{N+1})$. By \eqref{F1+}, we have
\begin{align}\label{F3+}
\nonumber
&2Nf^N-2(N+1)f^{N+1}\geq -2\langle g^{N+1}, (N+1)x^{N+1}-Nx^N- x^\star\rangle=\\
\nonumber   &\left\|  \sqrt[4]{\tfrac{B^2(N+1)}{R^2}}\left((N+1)x^{N+1}-Nx^N- x^\star\right)-\sqrt[4]{\tfrac{R^2}{B^2(N+1)}}g^{N+1} \right\|^2
\\
& -\tfrac{R}{B\sqrt{N+1}}\left\| g^{N+1} \right\|^2-\tfrac{B\sqrt{N+1}}{R} \left\|  (N+1)x^{N+1}-Nx^N- x^\star \right\|^2.
\end{align}
By multiplying \eqref{F2+} and \eqref{F3+} by $\tfrac{1}{2RB\sqrt{N+1}}$ and $\tfrac{1}{2(N+1)}$, respectively, and summing them, we get
\begin{align*}
f^{N+1}+\tfrac{N}{2BR\sqrt{N+1}}\left( f^{N}-\tfrac{BR}{\sqrt{N+1}} \right)^2&-\tfrac{NBR}{2(N+1)\sqrt{N+1}}\leq 
\tfrac{B}{2R\sqrt{N+1}}\left\| x^1-x^\star  \right\|^2
+
\\
& \tfrac{R}{2B(N+1)\sqrt{N+1}}\left\| g^{N+1} \right\|^2.
\end{align*}
As $\|g^{N+1}\|\leq B$ and $\|x^1-x^\star\|\leq R$, we obtain $f^{N+1}\leq \frac{BR}{\sqrt{N+1}}$, and the proof is complete. 
\end{proof}

\section{Convex feasibility problem}\label{Sec_alt_P}

A convex feasibility problem seeks to find a point in the intersection of convex sets. Let $C_1, C_2, \ldots, C_m$ be closed convex sets in $\mathbb{R}^n$. The convex feasibility problem can be formulated as finding $x \in \mathbb{R}^n$ such that
\begin{align}\label{F_P}
x \in \bigcap_{i=1}^{m} C_i.
\end{align}

We assume that $\bigcap_{i=1}^{m} C_i\neq \emptyset$. This problem arises in various applications, such as optimization, signal processing, and machine learning; see e.g. \cite{bauschke1996projection, combettes1996convex}. Let $d_C(.)$ denote the distance function to $C$, that is $d_C(x)=\min_{y\in C} \|x-y\|$. Moreover, $\Pi_C(.)$ stands for the projection on $C$. We introduce Algorithm \ref{F_OSM}, an adaptive greedy algorithm. Note that we call Algorithm \ref{F_OSM} the greedy algorithm when one is used as a step size instead of $1-\tfrac{k}{N+1}$.

\begin{algorithm}
\caption{Adaptive greedy algorithm}
\begin{algorithmic}
\smallskip
\State \textbf{Parameters:} number of iterations $N$
\smallskip
\State \textbf{Inputs:} convex sets $C_i$, $i\in\{1, \dots, m\}$, initial iterate $x^1 \in \mathbb{R}^n$ 
\smallskip
\State For $k=1, 2, \ldots, N$ perform the following steps:\\
\begin{enumerate}
\item
Select an $i_k\in \text{argmax}\ d_{C_{i}}(x^k)$.
\item
Set $x^{k+1}=x^k-\tfrac{N+1-k}{N+1}\left(x^k-\Pi_{C_{i_k}} (x^k)\right)$.
\end{enumerate}
\smallskip
\State \textbf{Output:} last iterate $x^{N+1}$ 
\end{algorithmic}
\label{F_OSM}
\end{algorithm}

It is seen problem \eqref{F_P} can be written as the following optimization problem 
\begin{align}\label{F_PP}
\min_{x\in\mathbb{R}^n} \left(\max_i d_{C_i}(x) \right).
\end{align}
Note that the objective function of problem \eqref{F_PP} is 1-Lipschitz continuous convex function. It is seen that $x=\Pi_C(y)$ implies that $\tfrac{1}{\|y-x\|}(y-x)\in\partial d_C(y)$. Hence, Algorithm \ref{F_OSM} is equivalent to  Algorithm \ref{OSM} applied to problem \eqref{F_PP}. Therefore, we have the following convergence rate.  

\begin{proposition}\label{T_OSM_F}
Let $C_1, C_2, \ldots, C_m$ be closed convex sets in $\mathbb{R}^n$. Consider $N$ iterations of Algorithm \ref{F_OSM} starting from an initial iterate $x^1 \in X$ satisfying $\|x^1 - x^\star\| \le R$ for some $x^\star\in \bigcap_{i=1}^{m} C_i$. 
We have 
\begin{align}
\max_i d_{C_i}(x^{N+1})\leq \frac{R}{\sqrt{N+1}}.
\end{align}
\end{proposition}

In a similar vein, Algorithm \ref{OSM+M} can be tailored for convex feasibility problems; see Algorithm \ref{F_OSM+M}. In this algorithm, the step sizes remain fixed regardless of the number of iterations, $N$. By using Theorem \ref{OSM+M}, we can infer the following proposition.    

\begin{algorithm}
\caption{Adaptive greedy algorithm with a momentum term}
\begin{algorithmic}
\smallskip
\State \textbf{Parameters:} number of iterations $N$
\smallskip
\State \textbf{Inputs:} convex sets $C_i$, $i\in\{1, \dots, m\}$, initial iterate $x^1 \in \mathbb{R}^n$ 
\smallskip
\State Set $x^0=x^1$. For $k=1, 2, \ldots, N$ perform the following steps:\\
\begin{enumerate}
\item
Select an $i_k\in \text{argmax}\ d_{C_{i}}(x^k)$.
\item
Set $x^{k+1}=x^k-\tfrac{1}{k+1}\left(x^k-\Pi_{C_{i_k}} (x^k)\right)+\tfrac{k-1}{k+1}\left(x^k-x^{k-1}\right)$.
\end{enumerate}
\smallskip
\State \textbf{Output:} last iterate $x^{N+1}$ 
\end{algorithmic}
\label{F_OSM+M}
\end{algorithm}

\begin{proposition}
Let $C_1, \dots, C_m$ be closed convex sets in $\mathbb{R}^n$. Consider $N$ iterations of Algorithm \ref{F_OSM+M} starting from an initial iterate $x^1 \in X$ satisfying $\|x^1 - x^\star\| \le R$ for some $x^\star\in \bigcap_{i=1}^{m} C_i$. 
We have 
\begin{align}
\max_i d_{C_i}(x^{N+1})\leq \frac{R}{\sqrt{N+1}}.
\end{align}
\end{proposition}

In what follows, we show that both Algorithm \ref{F_OSM} and \ref{F_OSM+M} attain the best convergence rate for convex feasibility problems with $m\geq N+1$. To this end, we only consider projection based methods that generate a sequence $\{x^k\}$ that satisfies the following assumption:\\
\textbf{Assumption A} The sequence $\{x^k\}$ satisfies 
\begin{align}
    x^{k+1}\in x^1+\text{Lin}\{x^1-\Pi_{C_{i_1}}(x^1), \dots, x^k-\Pi_{C_{i_k}}(x^k) \},
\end{align}
where $i_1, \dots, i_k\in\{1, \dots, m\}$.

The following theorem gives a resisting oracle in this setting. 

\begin{theorem}
Let $R>0$ and $N\in\mathbb{N}$ be given. There exist closed convex sets $C_1, \dots, C_{N+1}\subseteq \mathbb{R}^{N+1}$ and $x^1\in \mathbb{R}^{N+1}$ such that for any sequence $\{x^k\}$ satisfying Assumption A, we have 
$$
\min_{1\leq k\leq N+1}\left( \max_i d_{C_i}(x^k)\right)\geq \tfrac{R}{\sqrt{N+1}},
$$
where $R=d_{\cap_i C_i}(x^1)$.
\end{theorem}
\begin{proof}
Suppose that $C_i=\{x: \langle e_i, x\rangle=\tfrac{R}{\sqrt{N+1}}\}$ for $i\in \{1, \dots, N+1\}$  and $x^1=0$. It is readily seen that 
$$
\bigcap_{i=1}^{N+1} C_i=\{\tfrac{R}{\sqrt{N+1}}e\}, \ \ \ d_{\cap_i C_i}(x^1)=R.
$$
In addition, we have 
\begin{align}\label{L_B_F}
 \max_i d_{C_i}(x)\geq \tfrac{R}{\sqrt{N+1}}   
\end{align}
 for $x$ with at least one zero component. On the other hand, for $x\notin C_i$ with some zero components, we have 
$$
\supp( x)\subseteq \supp\left( x-\Pi_{C_i}(x)\right), \ \ 
\left| \supp\left( x-\Pi_{C_i}(x)\right) \right|\leq \left| \supp(x) \right|+1
$$
Note that if $x\in C_i$, we get $ x-\Pi_{C_i}(x)=0$. Hence, at least $N+2-k$ components of $x^k$ must be zero.  Inequality \eqref{L_B_F} implies 
$$
\max_i d_{C_i}(x^k)\geq \tfrac{R}{\sqrt{N+1}}, \ \ k\in\{1, \dots, N+1\},
$$
which completes the proof. 
\end{proof}

one can deduce that the convergence rate of the greedy algorithm in terms of the last iterate is $\mathcal{O}\left( \tfrac{1}{\sqrt[4]{N}}\right)$ by employing Theorem \ref{T_Polyak}. However, the given rate is very loose for problems with $m<N+1$ and it is advisable to directly investigate the method. To make our point clear, we establish that the convergence rate of the  alternating projection method is $\mathcal{O}\left( \tfrac{1}{\sqrt{N}}\right)$ in terms of the last iterate. Note that  the greedy algorithm amounts to the alternating projection when $m=2$ and Algorithm \ref{Alternating} describes the alternating projection method. 

\begin{algorithm}
\caption{Alternating projection method}
\begin{algorithmic}
\smallskip
\State \textbf{Parameters:} number of iterations $N$
\smallskip
\State \textbf{Inputs:} convex sets $C_1$ and $C_2$, initial iterate $x^1 \in C_2$ 
\smallskip
\State For $k=1, 2, \ldots, N$ perform the following step:\\
\begin{enumerate}
\item $x^{k+1}=\Pi_{C_2}\left( \Pi_{C_1}(x^k)\right)$
\end{enumerate}
\smallskip
\State \textbf{Output:} last iterate $x^{N+1}$ 
\end{algorithmic}
\label{Alternating}
\end{algorithm}

In the following theorem, we give a new convergence for the alternating projection in terms of the last iterate. 

\begin{theorem}\label{T_Alt}
Let $C_1, C_2\subseteq \mathbb{R}^n$ be closed and convex. Consider $N$ iterations of Algorithm \ref{Alternating}  with $\|x^1 - x^\star\| \le R$ for some $x^\star\in C_1\cap C_2$ and $x^1\in C_2$. Then 
\begin{align}\label{R_Alt}
d_{C_1}(x^{N+1})\leq R\sqrt{\frac{(2N)^{2N} }{ (2N+1)^{2N+1} }}.
\end{align}
\end{theorem}
\begin{proof}
Let $y^k=\Pi_{C_1} (x^k)$ for $k\in\{1, \dots, N+1\}$.  Without loss of generality, we assume that $x^\star=0$. As  $x^{k+1}=\Pi_{C_2} (y^k)$, by projection property, we have
\begin{align*}
\langle y^k-x^{k+1}, x^{k+1}-x^k\rangle\geq 0, \ \ \langle y^k-x^{k+1}, x^{k+1}\rangle\geq 0 \ \ k\in\{1, \dots, N\}.
\end{align*}
Multiplying these inequalities by $2k-1$ and $1+\tfrac{k}{N}$,respectively, and summing them leads to
\begin{align}\label{Alt_in1}
 \langle y^k-x^{k+1}, \tfrac{k(2N+1)}{N}x^{k+1}-(2k-1)x^k\rangle\geq 0.
\end{align}
On the other hand, we also have 
\begin{align}\label{Alt_in2}
\nonumber &\langle x^{k+1}-y^{k+1}, y^{k+1}-y^k\rangle\geq 0, \ & k\in\{1, \dots, N\},\\
&\langle x^k-y^k, y^k\rangle\geq 0 \ \ & k\in\{1, \dots, N+1\}.
\end{align}
Suppose that $s_k=\tfrac{2(2N)^{2(N-k)+1}}{(2N+1)^{2(N+1-k)}}$ for $k\in\{1, \dots, N+1\}$. By using \eqref{Alt_in1} and \eqref{Alt_in2}, we get
\begin{align*}
0\leq  &\sum_{k=1}^N s_k\langle y^k-x^{k+1}, \tfrac{k(2N+1)}{N}x^{k+1}-(2k-1)x^k\rangle+
\sum_{k=1}^N \tfrac{s_k(2N+2k-1)}{2N+1}\langle x^k-y^k, y^k\rangle+
\\
 & \sum_{k=1}^N \tfrac{4Nks_{k+1}}{2N+1}\langle x^{k+1}-y^{k+1}, y^{k+1}-y^k\rangle+
\tfrac{2}{2N+1}\langle x^{N+1}-y^{N+1}, y^{N+1}\rangle=\\
& \tfrac{(2N)^{2N} }{ (2N+1)^{2N+1} }\left\| x^1\right\|^2-\sum_{k=1}^{N} \tfrac{(2k-1)(2N)^{2(N+1-k)} }{ (2N+1)^{2(N+1-k)+1} }\left\| x^k-\tfrac{2N+1}{2N}x^{k+1}\right\|^2-\left\| x^{N+1}-y^{N+1}\right\|^2
\\
&-\sum_{k=1}^{N} ks_k\left\| y^k-\tfrac{2N+1}{2N}y^{k+1}\right\|^2.
\end{align*}
By using the above inequality, we get $\left\| x^{N+1}-y^{N+1}\right\|^2\leq \tfrac{(2N)^{2N} }{ (2N+1)^{2N+1} }\left\| x^1\right\|^2$ and the proof is complete. 
\end{proof}

It is worth mentioning that  the sequence $\{ \left(\tfrac{k}{k+1}\right)^k \}$ is a decreasing and tends to $\tfrac{1}{e}$. Thus, under the assumptions of Theorem \ref{T_Alt}, we get
\begin{align}\label{RR_Alt}
d_{C_1}(x^{N+1})\leq \tfrac{4R}{9\sqrt{2N+1} }.
\end{align}
Regarding the convergence rate of the alternating projection method, we have the following bound
\begin{align}\label{Amir_b}
\min_{1\leq k\leq N+1} d_{C_1}(x^k)\leq \tfrac{R}{\sqrt{N+1} },
\end{align}
see \cite[Corollary 8.22]{beck2017first}. The rate given in Theorem \ref{T_Alt} not only is tighter than \eqref{Amir_b}, but also it is given in  terms of the last iterate. 

It is worth noting that we obtained the proof of Theorem \ref{T_Alt} by using Performance estimation developed in \cite{taylor2017exact}. In the following theorem, we construct an example for which bound \eqref{R_Alt} is tight. 

\begin{theorem}\label{exact_Alt_e}
Let $N\in\mathbb{N}$ be given. There exist closed convex sets $C_1, C_2$ with $x^\star\in C_1\cap C_2$  for which Algorithm \ref{Alternating} with initial point $x^1\in C_2$ generates $x^{N+1}$ such that  
\begin{align}
d_{C_1}(x^{N+1})= \|x^1-x^\star\|\sqrt{\frac{(2N)^{2N} }{ (2N+1)^{2N+1} }}.
\end{align}
\end{theorem}
\begin{proof}
Consider two subspaces $C_1=\{x\in\mathbb{R}^2: x_2=\tfrac{1}{\sqrt{2N}}x_1\}$ and $C_2=\{x\in\mathbb{R}^2: x_2=0\}$. It is seen that $C_1\cap C_2=\{0\}$. In addition, 
$$
\Pi_{C_1} (x)=x_1
\begin{pmatrix}
    \tfrac{2N}{2N+1} & \tfrac{\sqrt{2N}}{2N+1}
\end{pmatrix}^T  \ \ x\in C_2, \ \ \ 
\Pi_{C_2} (x)=x_1. 
$$
Suppose that 
$x^1=\begin{pmatrix}
    1 & 0
\end{pmatrix}^T$.
Algorithm \ref{Alternating} with the initial point $x^1$ after $N$ iterations generates 
$x^{N+1}=\begin{pmatrix}
    \left(\frac{2N}{2N+1}\right)^N  & 0
\end{pmatrix}^T$. Moreover, we have $d_{C_1}(x^{N+1})=\sqrt{\frac{(2N)^{2N} }{ (2N+1)^{2N+1} }}$, which completes the proof. 
\end{proof}
For the example given in Theorem \ref{exact_Alt_e}, we have \( d_{C_1}(x^{k+1}) \leq d_{C_1}(x^k) \), \( k \in \{1, \dots, N\} \). This indicates that  bound  \eqref{R_Alt} is also tight in terms of the best iterate.


\bibliographystyle{siamplain}
\bibliography{references}

\begin{thebibliography}{10}

\bibitem{barre2021worst}
{\sc M.~Barr{\'e}}, {\em Worst-case analysis of efficient first-order methods},
  PhD thesis, Universit{\'e} Paris sciences et lettres, 2021.

\bibitem{barre2020complexity}
{\sc M.~Barr{\'e}, A.~Taylor, and A.~d’Aspremont}, {\em Complexity guarantees
  for polyak steps with momentum}, in Conference on learning theory, PMLR,
  2020, pp.~452--478.

\bibitem{bauschke1996projection}
{\sc H.~H. Bauschke and J.~M. Borwein}, {\em On projection algorithms for
  solving convex feasibility problems}, SIAM review, 38 (1996), pp.~367--426.

\bibitem{beck2017first}
{\sc A.~Beck}, {\em First-order methods in optimization}, SIAM, 2017.

\bibitem{berman1994nonnegative}
{\sc A.~Berman and R.~J. Plemmons}, {\em Nonnegative matrices in the
  mathematical sciences}, SIAM, 1994.

\bibitem{bertsekas2015convex}
{\sc D.~Bertsekas}, {\em Convex optimization algorithms}, Athena Scientific,
  2015.

\bibitem{boyd}
{\sc S.~Boyd, L.~Xiao, and A.~Mutapcic}, {\em Subgradient methods}, lecture
  notes of EE392o, Stanford University, Autumn Quarter, 2004 (2003),
  pp.~2004--2005.

\bibitem{cai2024last}
{\sc X.~Cai and J.~Diakonikolas}, {\em Last iterate convergence of incremental
  methods and applications in continual learning}, arXiv preprint
  arXiv:2403.06873,  (2024).

\bibitem{combettes1996convex}
{\sc P.~L. Combettes}, {\em The convex feasibility problem in image recovery},
  in Advances in imaging and electron physics, vol.~95, Elsevier, 1996,
  pp.~155--270.

\bibitem{davis2019stochastic}
{\sc D.~Davis and D.~Drusvyatskiy}, {\em Stochastic model-based minimization of
  weakly convex functions}, SIAM Journal on Optimization, 29 (2019),
  pp.~207--239.

\bibitem{defazio2023and}
{\sc A.~Defazio, A.~Cutkosky, H.~Mehta, and K.~Mishchenko}, {\em When, why and
  how much? adaptive learning rate scheduling by refinement}, arXiv preprint
  arXiv:2310.07831,  (2023).

\bibitem{defazio2024road}
{\sc A.~Defazio, H.~Mehta, K.~Mishchenko, A.~Khaled, A.~Cutkosky, et~al.}, {\em
  The road less scheduled}, arXiv preprint arXiv:2405.15682,  (2024).

\bibitem{drori2016optimal}
{\sc Y.~Drori and M.~Teboulle}, {\em An optimal variant of kelley's
  cutting-plane method}, Mathematical Programming, 160 (2016), pp.~321--351.

\bibitem{grimmer2023some}
{\sc B.~Grimmer and D.~Li}, {\em Some primal-dual theory for subgradient
  methods for strongly convex optimization}, arXiv preprint arXiv:2305.17323,
  (2023).

\bibitem{hazan2019revisiting}
{\sc E.~Hazan and S.~Kakade}, {\em Revisiting the polyak step size}, arXiv
  preprint arXiv:1905.00313,  (2019).

\bibitem{huang2024analytic}
{\sc Y.-K. Huang and H.-D. Qi}, {\em Analytic analysis of the worst-case
  complexity of the gradient method with exact line search and the polyak
  stepsize}, arXiv preprint arXiv:2407.04914,  (2024).

\bibitem{jain2019making}
{\sc P.~Jain, D.~M. Nagaraj, and P.~Netrapalli}, {\em Making the last iterate
  of {SGD} information theoretically optimal}, SIAM Journal on Optimization, 31
  (2021), pp.~1108--1130.

\bibitem{lang2012first}
{\sc S.~Lang}, {\em A first course in calculus}, Springer Science \& Business
  Media, 2012.

\bibitem{liu2023revisiting}
{\sc Z.~Liu and Z.~Zhou}, {\em Revisiting the last-iterate convergence of
  stochastic gradient methods}, arXiv preprint arXiv:2312.08531,  (2023).

\bibitem{liu2024last}
{\sc Z.~Liu and Z.~Zhou}, {\em On the last-iterate convergence of shuffling
  gradient methods}, arXiv preprint arXiv:2403.07723,  (2024).

\bibitem{last_iterate_sgd}
{\sc F.~Orabona}, {\em Last iterate of sgd converges even in unbounded
  domains}, 2020,
  \url{https://parameterfree.com/2020/08/07/last-iterate-of-sgd-converges-even-in-unbounded-domains/}.
\newblock Accessed: 2024-06-21.

\bibitem{polyak1969minimization}
{\sc B.~T. Polyak}, {\em Minimization of unsmooth functionals}, USSR
  Computational Mathematics and Mathematical Physics, 9 (1969), pp.~14--29.

\bibitem{taylor2017exact}
{\sc A.~B. Taylor, J.~M. Hendrickx, and F.~Glineur}, {\em Exact worst-case
  performance of first-order methods for composite convex optimization}, SIAM
  Journal on Optimization, 27 (2017), pp.~1283--1313.

\bibitem{zamani2023exact}
{\sc M.~Zamani and F.~Glineur}, {\em Exact convergence rate of the last iterate
  in subgradient methods}, arXiv preprint arXiv:2307.11134,  (2023).

\end{thebibliography}
\end{document}